\documentclass[aop,seceqn,MSNbibl,citesort,dvips]{arximspdf}

\doi{10.1214/10-AOP594}
\volume{39}
\issue{2}
\pubyear{2011}
\firstpage{417}
\lastpage{428}

\makeatletter

\newtheorem{theorem}{Theorem}[section]
\newtheorem{lemma}{Lemma}[section]
\newtheorem{proposition}{Proposition}[section]
\newtheorem{corollary}{Corollary}[section]

\newproclaim{Definition}{Definition}

\newcommand{\R}{\mathbb{R}}
\newcommand{\E}{\mathbb{E}}
\newcommand{\PP}{\mathbb{P}}
\newcommand{\C}{\mathbb{S}}
\newcommand{\Z}{\mathbb{Z}}

\def\bptnote#1{}

\makeatother

\begin{document}
\begin{frontmatter}

\title{T. E. Harris's contributions to recurrent Markov processes and
stochastic flows}
\runtitle{Recurrent Markov processes and stochastic flows}

\begin{aug}
\author[A]{\fnms{Peter} \snm{Baxendale}\corref{}\ead[label=e1]{baxendal@usc.edu}}
\runauthor{P. Baxendale}
\affiliation{University of Southern California}
\address[A]{Department of Mathematics\\
University of Southern California\\
3620 S. Vermont Avenue, KAP 108\\
Los Angeles, California 90089-2532\\
USA\\
\printead{e1}} 
\end{aug}

\received{\smonth{8} \syear{2010}}

%
\begin{abstract}
This is a brief survey of T. E. Harris's work on recurrent Markov
processes and on stochastic flows, and of some more recent work in
these fields.
\end{abstract}

%

\end{frontmatter}

\section{Introduction} I was a colleague of Ted Harris at USC, first
as a sabbatical visitor for the academic year 1982--1983 and later
as a regular faculty member from 1988 until his retirement. During
this time, I spent many hours in discussions with Ted, and
appreciated his insight into all areas of probability theory. I
feel deeply grateful for the opportunity to have learned so much
from him.

This paper covers two areas of Harris's work. Early in his career
he wrote his seminal paper \cite{Har-stat} on the existence and
uniqueness of stationary measures for Markov processes satisfying
a certain recurrence condition. Nowadays, this is called Harris
recurrence, although Ted was far too modest a person to ever use
the term himself. Section \ref{sec-rec} contains a brief account
of \cite{Har-stat} together with some indication of later
developments based on his idea.

Later in his career Harris became interested in stochastic flows.
These can be regarded as random mappings of an entire state space
into itself. Section \ref{sec-flow} starts with some details of
Harris's construction of a stochastic flow as a limit of random
stirring processes. (Here, we see a transition from Poisson point
processes and percolation theory into processes more frequently
described by stochastic differential equations.) The rest of this
section describes Harris's work on isotropic stochastic flows and
coalescing stochastic flows, together with some more recent
developments in these areas.
\eject

\section{Harris recurrence} \label{sec-rec}

Consider a (time homogeneous) discrete time Markov process $\{X_n\dvtx
n \ge0\}$ taking values in a measurable space $(S,\mathcal{B})$. For
simplicity here, we will assume $\mathcal{B}$ is separable, although
many results are valid in a more general setting. Let $P(x,
\cdot)$ denote the transition probability function, and let
$\PP^x$ denote the law of the process with initial condition $X_0=
x$. A nonzero measure $\mu$ on $(S,\mathcal{B})$ is a
\textit{stationary measure} for the Markov process if
\[
\mu(A) = \int_S \PP(x,A)\mu(dx) \qquad\mbox{for all } A \in
\mathcal{B}.
\]
\begin{Definition*}
For any set $A \in\mathcal{B}$, define the
hitting time $T_A = \inf\{n \ge1\dvtx X_n \in A\}$. Let $m$ be a
$\sigma$-finite measure on $(S, \mathcal{B})$. The Markov process
is said to be \textit{$m$-recurrent} if $\PP^x(T_A < \infty) = 1$ for
all $x \in S$ whenever $m(A)> 0$.
\end{Definition*}

Nowadays, we say the Markov process is \textit{Harris recurrent} if it
is $m$-recurrent for some nonzero $m$.
\begin{theorem}[(Harris \cite{Har-stat})] Assume the Markov
process $\{X_n\dvtx n \ge0\}$ is $m$-recurrent for some nonzero
$\sigma$-finite measure $m$. Then there is a stationary measure
$\mu$ for the Markov process. Moreover, $\mu$ is unique up to a
constant multiplier, and $m$ is absolutely continuous with respect
to $\mu$.
\end{theorem}

In the case of a discrete state space $S$, this result was already
known; see Derman \cite{Der} and Chung \cite{Chu}. The usual
definition of recurrence to a point $a \in S$ corresponds to
$m$-recurrence with $m$ taken to be the Dirac measure $\delta_a$ at
the point $a$, and then the measure
\[
\mu(B) = \sum_{n=1}^\infty\PP^a\bigl(X_n \in B, T_{\{a\}} \ge
n\bigr),\qquad
B \in\mathcal{B},
\]
is stationary. In this setting, the times of visits of $X_n$ to the
recurrent point $a$ are renewal times for the Markov process.
However, with a nondiscrete state space $S$ it will typically be
necessary to look for recurrence to larger sets $A$ than just
singleton sets, and then the locations of the visits within $A$ are
of interest.

For fixed $A \in\mathcal{B}$ with $m(A) > 0$, denote by $T_A^{(k)}$
the times of visits to $A$, so that $T_A^{(1)} = T_A$ and
$T_A^{(k+1)} = \inf\{n > T_A^{(k)}\dvtx X_n \in A\}$. The Markov
process $Z_m = X_{T_A^{(m)}}$ with values in $A$ is called the
\textit{process on $A$}. The transition probabilities for the
process on $A$ will be denoted $\PP_A(x,\cdot)$. The following
result, combining several lemmas in Harris \cite{Har-stat}, gives
a characterization of the stationary measure showing how the
discrete space result is extended.
\begin{proposition} \label{prop-mu} Assume the Markov process
$\{X_n\dvtx
n \ge0\}$ is $m$-recurrent for some nonzero
$\sigma$-finite measure $m$, and that $m(A) > 0$. Suppose that
$\mu_A$ is a stationary probability measure for the process on $A$
and that $m$ is absolutely continuous with respect to $\mu_A$ on
$A$. Define $\mu$ on $(S, \mathcal{B})$ by
%
%
\begin{equation} \label{mu}
\mu(B) = \int_{A}\Biggl( \E^x \sum_{n=1}^{T_A} 1_B(X_n) \Biggr)
\mu_A(dx),\qquad B \in\mathcal{B}.
\end{equation}
Then $\mu$ is a stationary measure for the Markov process $\{X_n\dvtx
n \ge0\}$ on $S$ and $m$ is absolutely continuous with respect to
$\mu$. Moreover, any other stationary measure is a multiple of
$\mu$.
\end{proposition}

It remains to consider the existence of the stationary probability
$\mu_A$ for the process on $A$. At this point, we see the
conflicting requirements for the set $A$: typically it has to be
larger than a singleton so that $m(A) > 0$, but it should be
chosen small enough so that the process on $A$ has some good
behavior.

Let $p^n(x,y)$ denote the absolutely continuous part of the $n$-step
transition probability $P^n(x,\cdot)$ with respect to $m$. The
following technical lemma of Harris is based in the idea that
$m$-recurrence implies that for each $x$, the set
$T(x) = \{y \in S\dvtx p^n(x,y) = 0 \mbox{ for all } n \ge1\}
$
must satisfy $m(T(x)) = 0$.
\begin{lemma} \label{lem-har} Assume the Markov process $\{X_n\dvtx n
\ge
0\}$ is $m$-recurrent for some nonzero
$\sigma$-finite measure $m$. For any $r \in(0,1)$, there exist $A
\in\mathcal{B}$ with \mbox{$0 < m(A) < \infty$}, a~positive integer $k$ and
a positive constant $s$ such that for all $x \in A$,
\[
m\{y \in A\dvtx p^1(x,y)+\cdots+ p^k(x,y) > s\} > r m(A).
\]
\end{lemma}

This enables Harris to obtain a Doeblin-like condition on the
transition operator $R(x,\cdot) = (P_A(x,\cdot)+ \cdots+
P_A^k(x,\cdot))/k$ and the existence of the stationary probability
$\mu_A$ follows directly.

\subsection{Small sets}

\begin{Definition*} A set $A \in\mathcal{B}$ is a \textit{small set}
if there exist a positive integer $k$, a probability measure $\nu$
and a constant $\beta> 0$ such that
%
%
\begin{equation} \label{small}
\PP^k(x, \cdot) \ge\beta\nu(\cdot) \qquad\mbox{for all }
x \in A.
\end{equation}
\end{Definition*}

The following result is a strengthening of Harris's Lemma
\ref{lem-har}.
\begin{proposition}[(Orey \cite{Ore})] \label{prop-orey} Assume
the Markov process $\{X_n\dvtx n \ge0\}$ is $m$-recurrent for some nonzero
$\sigma$-finite measure $m$. Every set $E \in\mathcal{B}$ such that
\mbox{$m(E) > 0$} contains a set $A \in\mathcal{B}$ such that $0 < m(A) <
\infty$ and
\[
\inf\{ p^k(x,y)\dvtx x,y \in A\} > 0
\]
for some positive integer $k$.
\end{proposition}
\begin{corollary} \label{cor-orey} The Markov process $\{X_n\dvtx n \ge
0\}
$ is $m$-recurrent for some nonzero
$\sigma$-finite measure $m$ if and only if there exists a small
set $A \in\mathcal{B}$ such that $\PP^x(T_A < \infty) = 1$ for all
$x \in S$.
\end{corollary}
\begin{pf} Any set $A$ with the property described in Proposition
\ref{prop-orey} is a small set, with $\nu(B) = m(A \cap B)/m(A)$.
Conversely, suppose $A$ has the properties described in Corollary
\ref{cor-orey}, with $k$, $\nu$ and $\beta$ as in (\ref{small}).
If $\nu(B) > 0$, then $\PP^x(X_{n+k} \in B) \ge\beta\nu(B)
\PP^x(X_n \in A)$. Since $\PP^x(T_A < \infty) = 1$ for all $x \in
S$, then $\PP^x(X_n \in A \mbox{ infinitely often}) = 1$ and so
$\PP^x(X_n \in B \mbox{ infinitely often}) = 1$. It follows that
$\{X_n\dvtx n \ge0\}$ is $\nu$-recurrent.
\end{pf}

This result shows that in some sense the study of Harris
recurrence is equivalent to the study of small sets with almost
surely finite hitting times. The property that $A$ is small can be
used in two distinct but related ways. For convenience assume
here that $A$ is a small set with $k = 1$. The general case can be
handled by applying the methods described below to the $k$-step
process $\{X_{nk}\dvtx n \ge0\}$.

\subsubsection{Coupling}
Two copies $\{X_n\dvtx n \ge0\}$ and $\{X'_n\dvtx n \ge0\}$ of the
Markov process can be coupled so that $\PP(X_{n+1} =
X'_{n+1}|X_n,X'_n) \ge\beta$ whenever\break $(X_n,X'_n) \in A \times
A$. This implies that two copies $\{Z_m\dvtx m \ge0\}$ and $\{Z'_m\dvtx
m \ge0\}$ of the process on $A$ can be coupled so that
$\PP(Z_{m+1} =Z'_{m+1} | Z_m, Z'_m) \ge\beta$. It follows easily
that the process on $A$ has a unique stationary probability
$\mu_A$, say, with
\[
\|\PP_A^m(x,\cdot) - \mu_A\| \le2(1-\beta)^m
\]
for all $m \ge0$ and $x \in A$. Details of this coupling
argument may be found in Lindvall \cite{Lin}.

\subsubsection{The split chain}

A set $C$ is said to be an \textit{atom} for the Markov chain $\{
X_n\dvtx
n \ge0\}$ if $\PP(x,\cdot) = \PP(y,\cdot)$ for all $x,y \in C$.
If $C$ is an atom, then the times of visits to $C$ are renewal
times for the Markov chain. The excursions away from $C$ will be
independent and identically distributed, and many questions about
the large time behavior of the Markov chain may be resolved using
this fact.

A singleton set $C = \{a\}$ is an clearly an atom. Nummelin
\cite{Num1} showed how to use a small set $A$ to build a Markov
chain with an atom. [More generally, Nummelin assumes a
minorization condition $\PP(x,C) \ge h(x)\nu(C)$ with some
nontrivial nonnegative function $h$; here, we specialize to $h(x)
= \beta1_A(x)$.] Consider the \textit{split chain}
$\{(X_{n},Y_{n})\dvtx n \geq0 \}$ with state space $S \times\{0,1\}$
and transition probabilities given by
\begin{eqnarray*}
P\{Y_{n} =1 | \mathcal{F}_{n}^{X} \vee\mathcal{F}_{n-1}^{Y} \} & = &
\beta1_{C}(X_{n}), \\
P\{X_{n+1} \in A | \mathcal{F}_{n}^{X}
\vee\mathcal{F}_{n}^{Y} \} & = & \cases{
\nu(A), &\quad if $Y_n = 1$,\vspace*{2pt}\cr
\displaystyle {\frac{P(X_n,A) - \beta
1_C(X_n) \nu(A)}{1 - \beta1_C(X_n)}}, &\quad if $Y_n =
0.$}
\end{eqnarray*}
Here, $\mathcal{F}_{n}^{X} = \sigma\{ X_{r} \dvtx 0 \leq r \leq n \}$ and
$\mathcal{F}_{n}^{Y} = \sigma\{ Y_{r} \dvtx 0 \leq r \leq n \}$. Thus,
the split chain evolves as follows. Given $X_n$, choose $Y_n$ so
that $P(Y_n = 1) = \beta1_C(X_n)$. If $Y_n = 1$ then $X_{n+1}$
has distribution $\nu$, whereas if $Y_n = 0$, then $X_{n+1}$ has
distribution $(P(X_n, \cdot) - \beta1_C(X_n)\nu)/(1 - \beta
1_C(X_n))$. The split chain $\{(X_n,Y_n)\dvtx n \ge0\}$ is designed
so that it has an atom $S \times\{1\}$, and so that its first
component $\{X_{n}\dvtx n \geq0 \}$ is a copy of the original Markov
chain.

Thus, renewal theory can be used to describe the large time and
stationary behavior of the split chain $\{(X_n,Y_n)\dvtx n \ge0\}$,
and hence of its first component $\{X_n\dvtx n \ge0\}$.

Much more information about small sets, and more generally about
the ergodic theory of Harris recurrent Markov chains, may be found
in the books of Revuz~\cite{Rev}, Nummelin \cite{Num2} and Meyn
and Tweedie \cite{MT-book} and the references therein.

\subsection{Positive recurrence and rates of convergence}

A Harris recurrent Markov process $\{X_n\dvtx n \ge0\}$ is to be
\textit{positive Harris recurrent} if the stationary measure $\mu$ can be
normalized to be a probability measure on $(S,\mathcal{B})$.

Suppose that $A$ is a small set with $0 < m(A) < \infty$, and
recall that $\mu_A$ denotes the stationary probability measure for
the process on $A$. The stationary measure $\mu$ given by
(\ref{mu}) has total mass
\[
\mu(S) = \int_A \E^x(T_A) \mu_A(dx).
\]
Thus, the issue of positive recurrence depends on estimates of the
expected hitting times $\E^x(T_A)$ for $x \in A$.
\begin{proposition}[(Tweedie \cite{Twe})] \label{prop-LF1}
Assume $A$ is a small set.
Suppose there exist a measurable function $V \ge0$ and constants
$c >0, k$ such that $PV(x) \le V(x) -c$ for $x \notin A$ and
$PV(x) \le k$ for $x \in A$. Then $\sup_{x \in A}\E^x(T_A)
\le1 + k/c$ and so the Markov process is positive
Harris recurrent.
\end{proposition}

The proof is based on the fact that the first inequality assumed
for the Lyapunov--Foster function $V$ implies that when $X_0
\notin A$ the process $V(X_n)+cn$ stopped at time $T_A$ is a
supermartingale. With stronger assumptions on the function $V$,
together with conditions to ensure aperiodicity, results may be
obtained concerning the rate of convergence of $P^n(x,\cdot)$ to
the stationary probability measure $\pi$, say. Several such
conditions are given by Meyn and Tweedie \cite{MT-disc,MT-book}.
The following result includes also stronger conditions on the
small set $A$ so as to ensure aperiodicity.
\begin{proposition}[(Meyn and Tweedie \cite
{MT-disc})] \label{prop-LF2} Assume
the set $A \in\mathcal{B}$ satisfies $P(x,\cdot) \ge\beta
\nu(\cdot)$ for all $x \in A$, where $\beta> 0$ and $\nu(A) = 1$.
Assume also there exist a measurable function $V \ge1$ and
positive constants $\lambda< 1$ and $k$ such that $PV(x) \le
\lambda V(x)$ for $x \notin A$ and $PV(x) \le k$ for $x \in A$.
Then the Markov process has a unique stationary probability
measure $\pi$, say. Moreover, there exist positive constants
\mbox{$\gamma< 1$} and $M$ with the property that
%
%
\begin{equation} \label{Verg}
\biggl| \int_S f(y) P^n(x,dy) - \int f(y) \pi(dy)\biggr| \le M V(x)
\gamma^n
\end{equation}
for all $x \in S$ and $n \ge0$ whenever $f\dvtx S \to\R$ is a
measurable function such that $|f(y)| \le V(y)$ for all $y \in
S$.
\end{proposition}

More details about estimates of the form (\ref{Verg}) can be found
in Chapter 16: $V$-uniform ergodicity of \cite{MT-book}. In Meyn
and Tweedie \cite{MT-comp} it is shown that the constants $M$ and
$\gamma$ can be chosen depending only on $\lambda$, $k$ and
$\beta$. See also Baxendale \cite{Bax-comp}.

\subsection{Continuous time Markov processes}

\subsubsection{Sampled chains}

Suppose that $\{X_t\dvtx t \ge0\}$ is a time homogeneous Markov
process, and that $\{T(n)\dvtx n \ge0\}$ is an independent undelayed
renewal process with increment distribution $a$, for some
probability distribution on $(0,\infty)$. Then the sampled chain
$\{Y_n\dvtx n \ge0\}$ defined by $Y_n = X_{T(n)}$ is a time
homogeneous Markov chain. If $P_t(x,\cdot)$ denotes the time $t$
transition probability function for $\{X_t\dvtx t \ge0\}$, then
$\{Y_n\dvtx n \ge0\}$ has transition probability function
$K_a(x,\cdot) = \int P_t(x,\cdot)a(dt)$.

The $\Delta$-skeleton chain $Y_n = X_{n\Delta}$ corresponds to the
deterministic $a = \delta_\Delta$. Alternatively, the resolvent
chain observes the process $X$ at the times of a rate 1 Poisson
process and has transition probability
\[
R(x,\cdot) = \int P_t(x,\cdot)e^{-t}\,dt,
\]
which is the resolvent of the original continuous time process.
The advantages and disadvantages of various choices for $a$, and
the connections between the theories of Harris recurrence for
continuous time and discrete time Markov processes are discussed
in the papers \cite{MT-doeb,MT-cont} of Meyn and Tweedie. Results
analogous to those in Propositions \ref{prop-LF1} and
\ref{prop-LF2}, involving the action of the infinitesimal
generator $L$ on $V$, are given in Meyn and Tweedie \cite{MT-III}.

\subsubsection{Stopping times}

An alternative approach to recurrence for continuous time Markov
processes was developed independently by Maruyama and Tanaka
\cite{MaT} and Khas'minskii \cite{Kha}. Assume that $\{X_t\dvtx t
\ge0\}$ is a strong Markov process with right continuous paths
and left limits on a separable metric space $S$. Suppose that
$D_1$ and $D_2$ are open sets with disjoint closures
$\overline{D}_1$ and $\overline{D}_2$ with the property that the
stopping times $T_{D_1} = \inf\{t \ge0\dvtx X_t \in D_1\}$ and
$T_{D_2} = \inf\{t \ge0\dvtx X_t \in D_2\}$ are both $\PP^x$-almost
surely finite for all $x \in S$. Define inductively sequences
$\sigma_n$ and $\tau_n$ of stopping times by $\sigma_0 = \inf\{t
\ge0\dvtx X_t \in D_1\}$, $\tau_n = \inf\{t \ge\sigma_n\dvtx X_t \in
D_2\}$ for $n \ge0$, and $\sigma_n = \inf\{t \ge\tau_{n-1}\dvtx X_t
\in D_1\}$ for $n \ge1$. Then $\{Y_n \equiv X_{\sigma_n}\dvtx n \ge
0\}$ is a time homogeneous Markov chain on $\overline{D}_1$.

Assume for the moment that the process $\{Y_n\dvtx n \ge0\}$ on
$\overline{D}_1$ has a stationary probability measure
$\mu_{\overline{D}_1}$. Then it is shown in \cite{MaT} and
\cite{Kha} that the measure $\mu$ on $(S, \mathcal{B})$ defined by
%
%
\begin{equation} \label{mu2}
\mu(B) = \int_{\overline{D}_1}\biggl( \E^x \int_0^{\sigma_1}
1_B(X_s) \,ds \biggr) \mu_{\overline{D}_1}(dx),\qquad B \in
\mathcal{B},
\end{equation}
is a stationary measure for the original continuous time Markov
process $\{X_t\dvtx\break t \ge0\}$.

So far this is a very natural extension to continuous time of
Harris's result in Proposition \ref{prop-mu}. It remains to show
the existence of the stationary probability measure
$\mu_{\overline{D}_1}$, and this is where both \cite{MaT} and
\cite{Kha} impose extra conditions on the process $\{X_t\dvtx t \ge
0\}$. In particular, they both use additional properties of the
hitting distribution $\PP^x(X_{T_{D_1}} \in B)$ for $x \in
\overline{D}_2$ and $B \subset\overline{D}_1$ which ensure that
the process $\{Y_n\dvtx n \ge0\}$ satisfies Doeblin's condition (D).

The representation (\ref{mu2}) can be used to convert occupation
time estimates for the process $\{X_t\dvtx t \ge0\}$ in a very direct
way into estimates on the stationary measure~$\mu$, see, for
example, Baxendale and Stroock \cite{BS} and Baxendale~\cite{Bax-hopf}.

\section{Stirring processes and stochastic flows} \label{sec-flow}

\subsection{Random stirring in ${\R}^d$}

Let $\phi\dvtx \R^d \to\R^d$ be a homeomorphism of $\R^d$ onto itself
such that $\phi(x) = x$ whenever $\|x\| \ge K$, for some $K$. The
mapping $\phi$ can be thought of as a stirring of $\R^d$ centered
at the origin $0 \in\R^d$. (Harris used the term ``stirring''
originally in the case where the homeomorphism $\phi$ is volume
preserving, but it is convenient to keep the term in this more
general setting.) For $a \in\R^d$, the translated mapping
$\phi^a(x) = a+\phi(x-a)$ represents stirring centered at $a$.

Consider a Poisson point process on $\R^d \times(0,\infty)$ with
intensity $\lambda \,dx \,dt$. For each atom $(a,t)$ of the point
process, apply the mapping $\phi^a$ at time $t$. Then for $0 \le
s \le t < \infty$ the value $X_{st}$ of the stirring process is
the random mapping of $\R^d$ to itself obtained as the composition
of the stirrings $\phi^a$ at times $u$ for all the atoms $(a,u)$
with $s < u \le t$. A percolation argument, using the fact that
$\phi(x)-x$ has bounded support, can be used to show that for
sufficiently small $t-s$ the restriction of the mapping $X_{st}$
to any bounded set is almost surely given by the composition of a
finite number of stirrings. It follows that the process
$\{X_{st}\dvtx 0 \le s \le t < \infty$ is well defined.

The essence of this construction can be seen in the paper by Harris
\cite{Har-excl} on the construction of an exclusion process with
nearest neighbor rates. In \cite{Har-excl}, he considers a point
process on the set of bonds of the integer lattice $\Z^d$ and the
corresponding stirring switches the two ends of the bond. The
construction in \cite{Har-excl} allows the rates to depend on the
local configuration, but in the simplest case of constant rate it
fits into the setting above. Random stirring on the real line is
studied in the paper \cite{Lee} by Harris's student W. C. Lee.

\subsubsection{Convergence to a stochastic flow}

Consider the effect of letting the magnitude of the displacement
involved in each stirring $\phi$ tend to zero while letting the
intensity of the Poisson process tend to infinity. More
precisely, fix a compactly supported vector field $V$ on $\R^d$
and let $\phi_n$ denote the time $1/\sqrt{n}$ flow along $V$. Let
$\{X^n_{st}\dvtx 0 \le s \le t < \infty\}$ be the random stirring
process obtained using the stirring function $\phi_n$ together
with a Poisson process with rate $n\lambda \,dx \,dt$. Assume the
centering condition
\[
\int_{\R^d} V(x) \,dx = 0.
\]
Then under appropriate smoothness conditions on the vector field
$V$ the processes $\{X^n_{st}\dvtx 0 \le s \le t < \infty\}$ converge
weakly to a process $\{X_{st}\dvtx 0 \le s \le t < \infty\}$ with
values in the group of homeomorphisms of $\R^d$. This is proved in
Harris \cite{Har-hom} for the case $d=2$ when $V$ is divergence
free and rotationally symmetric. A more general form of result
(although on a compact manifold) is given in Matsumoto and
Shigekawa \cite{MaS}.

The process $\{X_{st}\dvtx 0 \le s \le t < \infty\}$ has the
properties:
\begin{enumerate}[(iii)]
\item[(i)] for each $x \in\R^d$ and $s \ge0$ the mapping $t \to
X_{st}(x)$ is continuous;
\item[(ii)] $X_{tu} \circ X_{st} = X_{su}$ whenever $s \le t \le
u$;
\item[(iii)] if $s_1 \le t_1 \le s_2 \le t_2 \le\cdots\le s_n \le
t_n$ the mappings $X_{s_it_i}$, $1 \le i \le n$ are independent;
\item[(iv)] the distribution of $X_{st}$ depends only on $t-s$.
\end{enumerate}
Any process with these properties will be called a
(time-homogeneous)
\textit{stochastic flow} on $\R^d$.
Much information about stochastic flows may be found in the books
of Kunita \cite{Kun} and Arnold \cite{Arn}.

The law of a stochastic flow is determined by the laws of its
$k$-point motions
\[
t \to(X_{0t}(x_1),X_{0t}(x_2),
\ldots,X_{0t}(x_k)) \in\R^{dk}
\]
for all $k \ge1$ and
$x_1,x_2,\ldots,x_k \in\R^d$. It is easy to see that each
$k$-point motion is a Markov process on $\R^{dk}$. Under suitable
regularity conditions, the $k$-point motions are diffusion processes
and the infinitesimal generator for the $k$ point motion can be
explicitly written in terms of the generator $L$, say, for the
one-point motion and a covariance matrix $B(x,y)$ for the two-point
motion. In particular for $f \in C^2(\R^d)$ with compact support,
\[
Lf(x) = \lim_{t \searrow0} \frac{\E f(X_{0t}(x))- f(x)}{t},\qquad
x \in\R^d,
\]
and
\[
B^{pq}(x,y) = \lim_{t \searrow0} \frac{\E[(X^p_{0t}(x)-
x^p)(X_{0t}^q(y)-y^q)]}{t},\qquad x,y \in\R^d.
\]
The operator $L$ and the matrix function $B$ are related by the
fact that $B(x,x)$ is the symbol of the operator $L$. Together,
$L$ and $B$ are called the local characteristics of the flow; see
Le Jan and Watanabe \cite{LeW}. For the stochastic flow
constructed above as the limit of random stirring processes, the
operator $L$ has constant coefficients and $B(x,y)$ depends only
$x-y$. This implies that the law of the stochastic flow is
homogeneous in space as well as time. For the example on $\R^2$
considered by Harris in \cite{Har-hom}, the rotational invariance
of the vector field $V$ implies that law of the stochastic flow is
invariant under rigid motions of $\R^2$, and in particular the
one-point motion is Brownian motion (up to a scaling factor).

\subsection{Isotropic stochastic flows}

A stochastic flow $\{X_{st}\dvtx 0 \le s \le t < \infty\}$ is
\textit{isotropic} if its law is invariant under rigid motions of $\R^d$.
Harris \cite{Har-hom} studied incompressible isotropic stochastic
flows on $\R^2$, and Baxendale and Harris \cite{BH} and Le Jan
\cite{LeJ} studied the general $d$-dimensional case.

For an isotropic stochastic flow, the generator $L$ for the
one-point motion is a multiple of the Laplace operator $\Delta$,
and the law of the flow is determined by the covariance matrix
$B$. Invariance under translations implies $B(x,y) = B(x-y,0)$,
and then invariance under rotations implies $B(x) \equiv B(x,0)$
satisfies $B(x) = G^*B(Gx)G$ for all real orthogonal matrices $G$.
This condition gives a representation of $B$ using Bessel
functions; see Yaglom \cite{Yag} and It\^o \cite{Ito}. A
corresponding representation for isotropic stochastic flows on a
sphere $\C^d$ appears in Raimond~\cite{Rai}.

The isotropy condition implies that certain geometric properties
of the flow can be calculated explicitly. For example, the length
$\|v_t\|$ of a tangent vector $v_t = DX_{0t}(x)(v)$ is a geometric
Brownian motion and the top Lyapunov exponent $\lambda_1 = \lim_{t
\to\infty} t^{-1} \log\|v_t\|$ can be calculated explicitly in
terms of $B$. Other local geometric properties such as the
curvature of a submanifold of $\R^d$ have been calculated; see Le
Jan \cite{LeJ-geom} and Cranston and Le Jan \cite{CrLeJ-geom}.

Of more interest are results involving the joint behavior of
infinitely many points. A result of Baxendale and Harris on the
length of a small curve in the case $\lambda_1 < 0$ has recently
been sharpened by Dimitroff \cite{Dim}. Results of Cranston,
Scheutzow and Steinsaltz \cite{CSS-iso,CSS-lin} show that, while
$\|X_{0t}(x)\|$ grows like $\sqrt{t}$ for each fixed $x$, if $D$
is a nonsingleton connected set in $R^d$ for $d \ge2$ and the
isotropic stochastic flow has $\lambda_1 > 0$ then $\sup\{
\|X_{0t}(y)\|\dvtx y \in D\}$ grows almost surely linearly as $t \to
\infty$.

For any measure (distribution of mass) $\nu$, let $\nu_t$ denote
the induced random measure $\nu\circ X_{0t}^{-1}$. Zirbel
\cite{Zir} contains estimates on the first two moments of $\nu_t$.
Recently Cranston and Le Jan \cite{CrLeJ-clt} and Dimitroff and
Scheutzow \cite{DS} have proved asymptotic normality of the
rescaled random measure $A \to\nu_t(\sqrt{t}A)$.

\subsection{Coalescing flows}

For vector fields $V_0, V_1, V_2,\ldots$ on $\R^d$ and
independent scalar Brownian motions $\{W_t^1\dvtx t \ge0\}, \{
W_t^2\dvtx
t \ge0\}, \ldots,$ consider the stochastic differential equation
%
%
\begin{equation} \label{sde}
dx_t = V_0(x_t) \,dt + \sum_{\alpha\ge1}
V_\alpha(x_t)\,dW_t^\alpha.
\end{equation}
Under suitable regularity and growth conditions on the vector
fields $V_0, V_1, V_2, \ldots,$ the strong solutions of
(\ref{sde}) for different initial conditions can be pieced
together to give a stochastic flow $\{X_{st}\dvtx 0 \le s \le t <
\infty\}$ of homeomorphisms $\R^d$; see, for example, Kunita
\cite{Kun}. The local characteristics of the flow are the operator
\[
Lf(x) = \sum_{p=1}^d V_0^p(x)\,\frac{\partial f}{\partial x^p}(x) +
\frac{1}{2} \sum_{p,q=1}^d \sum_{\alpha\ge1} V_\alpha
^p(x)V_\alpha^q(x)\,
\frac{\partial^2 f}{\partial x^p \,\partial x^q}(x)
\]
and the covariance function
%
%
\begin{equation} \label{B}
B^{pq}(x,y) = \sum_{\alpha\ge1} V_\alpha^p(x)V_\alpha^q(y).
\end{equation}
Conversely, any stochastic flow in which $L$ and $B$ have
sufficiently smooth coefficients arises as the solution of a
stochastic differential equation (taking the $V_\alpha$ to be an
orthonormal basis of the reproducing kernel Hilbert space of $B$)
and the flow consists of homeomorphisms.

Harris \cite{Har-coal} introduced the study of coalescing
stochastic flows. These are ones where the mappings $X_{st}$ may
be many to one. Harris studied the case $d=1$ with continuous
homogenous (in space) covariance function $B$, and obtained
conditions for coalescence in terms of the modulus of continuity
of $B$ at $0$. (In contrast, the ``Arratia flow'' of independent
coalescing Brownian motions; see \cite{Arr}, has discontinuous $B
= 1_{\{0\}}$.)

The issue of the existence of nonhomeomorphic stochastic flows in
dimensions $d \ge2$ was addressed by Darling \cite{Dar}. More
recently, Le Jan and Raimond \cite{LeJR-1,LeJR-2} have developed
new techniques to interpret the stochastic differential equation
(\ref{sde}) when the covariance function $B$ given by (\ref{B}) is
non-Lipschitz. In this more general setting, there is not only
the possibility of coalescence; there is also the possibility that
the solution of (\ref{sde}) has to be interpreted as a flow of
probability kernels. The flow of probability kernels, rather than
a flow of mappings, corresponds to the lack of uniqueness in the
solutions of (\ref{sde}). Examples of such flows include flows on
Euclidean space $\R^d$ and spheres $\C^d$ with isotropic, but
non-Lipschitz, covariance functions~$B$.

%

%
\printaddresses

\end{document}